\documentclass[10pt,a4paper,leqno]{amsart}
\usepackage[english]{babel}
\usepackage[T1]{fontenc}
\usepackage{amssymb}
\usepackage{amsfonts}
\usepackage{amsmath}
\usepackage[T1]{fontenc}
\usepackage{eucal}

\begin{document}
\title{\Large Domination conditions for  families of quasinearly subharmonic functions}
\author{Juhani Riihentaus}
\date{June 19, 2011}
\maketitle

\vspace{-0.35in}

\begin{center}
{Department of Mathematical Sciences,  University of Oulu\\
P.O. Box 3000, 90014 Oulun Yliopisto, Finland \\
juhani.riihentaus@gmail.com}
\end{center}

\vspace*{1ex}

\noindent{\emph{Abstract:}} Domar has  given  a condition that ensures the existence of  the largest subharmonic minorant of a given function. Later Rippon pointed out that a modification of Domar's argument gives in fact a better result. Using our previous, rather general and flexible, modification of Domar's original argument, we extend their results both to the subharmonic and quasinearly subharmonic settings.   

\vspace{0.5ex}

\section{Introduction}
\subsection{Results of Domar and Rippon} Suppose that $D$ is  a domain of ${\mathbb{R}}^n$, $n\geq 2$. Let $F:\, D\rightarrow [0,+\infty ]$ be an upper semicontinuous function. Let ${\mathcal{F}}$ be a family of subharmonic  functions $u: D\rightarrow [0,+\infty )$ which satisfy 
\[u(x)\leq F(x) \]
for all $x\in D$. 
Write 
\[ w(x):=\sup_{u\in {\mathcal{F}}}u(x), \,\, x\in D,\]
and let $w^*:D\rightarrow [0,+\infty ]$ be the upper semicontinuous regularization of $w$, that is 
\[ w^*(x):=\limsup_{y\rightarrow x}w(y).\]

Domar gave  the following result:
\subsection*{Theorem~A} (\cite{Do57},  Theorem~1 and Theorem~2,  pp.~430-431) \emph{If for some $\epsilon >0$,
\begin{equation}\int\limits_D [\log^+F(x)]^{n-1+\epsilon}\, dm_n(x)<+\infty ,
\end{equation}
then $w$ is locally bounded above in $D$, and thus $w^*$ is subharmonic in $D$.}

See \cite{Do57}, Theorem~1 and  Theorem~2,  pp.~430-431, see also \cite{Ni95}, pp.~67-69. As Domar points out, the original case of subharmonic functions in the result of \cite{Do57}, Theorem~1, p.~430,  is  due to Sj\"oberg \cite{Sj38} and Brelot \cite{Br45} (cf. also Green \cite{Gr52}). Observe, however, that Domar also sketches  a new proof for this Theorem~1 which uses  elementary methods and  applies   to  more general functions.  

Rippon generalized   Domar's  result in the following form:
\subsection*{Theorem~B} (\cite{Rip81},  Theorem~1,  p.~128) \emph{Let $\varphi :[0,+\infty ]\rightarrow [0,+\infty ]$ be an increasing function such that 
\begin{equation*}\int\limits_1^{+\infty} \frac{dt}{[\varphi (t)]^{1/(n-1)}}<+\infty .
\end{equation*}
If 
\begin{equation}\int\limits_D\varphi ( \log^+F(x))\, dm_n(x)<+\infty ,
\end{equation}
then $w$ is locally bounded above in $D$, and thus $w^*$ is subharmonic in $D$.}

As pointed out by Domar, \cite{Do57}, pp.~436-440,  by Nikolski, \cite{Ni95}, p.~69, and by Rippon, \cite{Rip81}, p.~129, the above results are for many particular cases sharp. For related results, see also \cite{Yo86}.

\vspace{1ex}

As Domar points out, \cite{Do57}, p. 430, the result of his Theorem~A  holds in fact for more general functions, that is, for functions which by good reasons might be -- and indeed already have been -- called quasinearly subharmonic functions. See 1.3 below for the definition of this function class. In addition, Domar has  given a related result for an even more  general function class $K(A,\alpha )$, where the above conditions (1) and (2) are replaced by a certain integrability condition on the decreasing rearrangement of $\log F$, see \cite{Do88}, Theorem~1, p.~485. Observe, however,  that in the case $\alpha =n$ Domar's class $K(A,n)$ equals with the class of nonnegative quasinearly subharmonic functions: If $u\in K(A,n)$, then $u$ is  $\nu_nA^{n+1}$-quasinearly subharmonic, and conversely, if $u\geq 0$ is $C$-quasinearly subharmonic, then $u\in K(C,n)$. 

\vspace{1ex}

Below we give a  general and at the same time flexible result which includes both Domar's and Rippon's results, Theorems~A and B above. See Theorem~2.1 and Corollary~2.5 below. For  previous preliminary, more or less standard results, see also \cite{He71}, Theorem~2 (d), p.~15, \cite{ArGa01}, Theorem~3.7.5, p.~83, \cite{Ri89}, Theorem~2, p.~71, and \cite{Ri07}, Theorem~2.2 (vi), p.~55 (see {\textbf{1.5}}~(v) below). 
\subsection{Notation} Our notation is rather standard, see e.g. \cite{Ri07} and \cite{He71}. For the convenience of the reader we, however, recall the following. $m_n$ is the Lebesgue measure 
in the Euclidean space ${\mathbb{R}}^n$, $n\geq 2$, and $\nu _n$ is  the Lebesgue measure of the unit ball $B^n(0,1)$ 
in ${\mathbb{R}}^n$, thus $\nu _n=m_n(B^n(0,1))$. $D$ is always a  domain in ${\mathbb{R}}^n$.
Constants will be denoted by $C$ and $K$. They are always nonnegative and may vary from line to line.
\subsection{Subharmonic functions and generalizations} We recall that an upper semicontinuous function $u:\, D\rightarrow [-\infty ,+\infty )$ is \emph{subharmonic} if 
for all closed balls $\overline{B^n(x,r)}\subset D$,
\[u(x)\leq \frac{1}{\nu _n\, r^n}\int\limits_{B^n(x,r)}u(y)\, dm_n(y).\]
The function $u\equiv -\infty $  is considered  subharmonic. 

We say that a function 
$u:\, D\rightarrow [-\infty ,+\infty )$ is \emph{nearly subharmonic}, if $u$ is Lebesgue measurable, $u^+\in {\mathcal{L}}^1_{\textrm{loc}}(D)$, 
and for all $\overline{B^n(x,r)}\subset D$,  
\begin{equation*}u(x)\leq \frac{1}{\nu _n\, r^n}\int\limits_{B^n(x,r)}u(y)\, dm_n(y).\end{equation*}
Observe that in the standard definition of nearly subharmonic functions one uses the slightly stronger  assumption that 
$u\in {\mathcal{L}}^1_{\textrm{loc}}(D)$, see e.g. \cite{He71}, p.~14. However, our above, slightly 
more general definition seems to be  more practical, see  e.g. \cite{Ri07},  Proposition~2.1~(iii) and Proposition~2.2~(vi) and (vii), pp.~54-55, and also {\textbf{1.5}} (i)--(vi) below. The following lemma emphasizes this fact still more:
\subsection{Lemma} (\cite{Ri07}, Lemma, p.~52) \emph{Let $u:\,D\rightarrow [-\infty ,+\infty )$ be Lebesgue measurable. 
Then  $u$ is nearly subharmonic (in the sense defined above) if and only if there exists a function $u^*$, subharmonic in $D$  such that $u^*\geq u$ and $u^*=u$ almost 
everywhere in $D$. Here $u^*$ is the upper semicontinuous regularization of $u$:
\begin{displaymath}u^*(x)=\limsup_{x'\rightarrow x}u(x').\end{displaymath}
}

\vspace{1ex} 
 
The proof follows at once  from \cite{He71}, proof of Theorem~1, pp.~14-15,  (and referring also to \cite{Ri07}, Proposition~2.1~(iii) and Proposition~2.2~(vii), pp.~54-55).
 
We say that a Lebesgue measurable function $u:\,D \rightarrow 
[-\infty ,+\infty )$ is \emph{$K$-quasinearly subharmonic}, if  $u^+\in{\mathcal{L}}^{1}_{\textrm{loc}}(D)$ and if there is a 
constant $K=K(n,u,D)\geq 1$
such that for all  $\overline{B^n(x,r)}\subset D$,    
\begin{equation} u_M(x)\leq \frac{K}{\nu _n\,r^n}\int\limits_{B^n(x,r)}u_M(y)\, dm_n(y)\end{equation}
for all $M\geq 0$, where $u_M:=\max\{u,-M\}+M$. A function $u:\, D\rightarrow [-\infty ,+\infty )$ is \emph{quasinearly subharmonic}, if $u$ is 
$K$-quasinearly subharmonic for some $K\geq 1$.

\vspace{1ex}

A Lebesgue  measurable function 
$u:\,D \rightarrow [-\infty ,+\infty )$ is \emph{$K$-quasinearly subharmonic n.s. (in the narrow sense)}, if $u^+\in{\mathcal{L}}^{1}_{\textrm{loc}}(D)$ and if there is a 
constant $K=K(n,u,D)\geq 1$
such that for all $\overline{B^n(x,r)}\subset D$,    
\begin{equation} u(x)\leq \frac{K}{\nu _n\,r^n}\int\limits_{B^n(x,r)}u(y)\, dm_n(y).\end{equation}
A function $u:\, D\rightarrow [-\infty ,+\infty )$ is \emph{quasinearly subharmonic n.s.}, if $u$ is 
$K$-quasinearly subharmonic n.s. for some $K\geq 1$. 

\vspace{1ex}

As already pointed out, Domar, \cite{Do57} and \cite{Do88},  considered nonnegative quasinearly subharmonic functions. Later on quasinearly subharmonic functions (perhaps with a different terminology, and sometimes in certain special cases, or just the corresponding generalized mean value inequality (3) or (4)) have been considered in many papers, see e.g. \cite{FeSt72}, Lemma~2, p.~172,  \cite{To86}, pp.~188-191, \cite{Ri89}, Lemma, p.~69, \cite{Pa88}, and \cite{Pa94}, Theorem~1, p.~19. For a rather detailed list of references, see \cite{Ri08$_2$}, and, for some more recent articles,  e.g. \cite{Ko07}, \cite{Ri07},  \cite{PaRi08}, \cite{Pa08}, \cite{Ri08$_1$},  \cite{Ri09},  \cite{PaRi10}, \cite{DoRi10$_1$}, \cite{DoRi10$_2$} and \cite{Ri10}. 

We recall here only that this  function class 
includes, among \mbox{others,}  subharmonic functions, and, more generally,  quasisubharmonic (see e.g. \cite{Br38},  \cite{Le45}, p.~309, \cite{Av61}, p.~136, \cite{He71}, p.~26)
 and also 
nearly subharmonic functions (see e.g. \cite{Br65}, pp.~30-31, \cite{He71}, p.~14),    also functions satisfying certain natural  growth conditions, especially  
certain eigenfunctions, and  polyharmonic functions. Also, the class of Harnack functions is included, thus, among others, nonnegative harmonic functions 
as well as nonnegative solutions of some elliptic equations. In particular, the partial differential equations associated with quasiregular mappings 
belong to this family of elliptic equations, see \cite{Vu82}.  
\subsection{}
For the sake of 
convenience of the reader we recall the following, see \cite{Ri07}, Proposition~2.1 and Proposition~2.2, pp.~54-55:  
\begin{itemize}
\item[(i)] \emph{A $K$-quasinearly subharmonic function n.s. is $K$-quasinearly subharmonic, but not necessarily conversely.}
\item[(ii)] \emph{A nonnegative Lebesgue measurable function is $K$-quasinearly subharmonic if and only if it is $K$-quasinearly subharmonic n.s.} 
\item[(iii)] \emph{A Lebesgue measurable function is \mbox{$1$-quasinearly} subharmonic if and only if it is \mbox{$1$-quasinearly} subharmonic n.s. and 
if and only if it is nearly subharmonic (in the sense  defined above).}
\item[(iv)] \emph{If $u:\, D\rightarrow [-\infty ,+\infty )$   is \mbox{$K_1$-quasinearly} subharmonic and 
$v:\, D\rightarrow [-\infty ,+\infty )$ is \mbox{$K_2$-quasinearly} subharmonic, then $\max \{u,v\}$ is $\max \{K_1,K_2\}$-quasinearly 
subharmonic in $D$. Especially,  
 $u^+:=\max\{u,0\}$ is $K_1$-quasinearly subharmonic in $D$.}
\item[(v)] \emph{Let ${\mathcal{F}}$ be a family of  $K$-quasinearly subharmonic (resp. $K$-quasinearly subharmonic n.s.) functions in $D$ and let 
$w:=\sup_{u\in {\mathcal{F}}}u$. If $w$ is Lebesgue measurable and} $w^+\in {\mathcal{L}}_{{\textrm{loc}}}^1(D)$, \emph{then $w$ is 
$K$-quasinearly subharmonic (resp. $K$-quasinearly subharmonic n.s.) in $D$.}
\item[(vi)] \emph{If  $u:\, D\rightarrow [-\infty ,+\infty )$   is quasinearly subharmonic n.s., then either $u\equiv -\infty $ or $u$ is finite almost 
everywhere in $D$, and} 
$u\in {\mathcal{L}}^1_{\textrm{loc}}(D)$.
\end{itemize}
\section{The result} 
\subsection{Theorem}  \emph{Let $K\geq 1$.  Let $\varphi :[0,+\infty ]\rightarrow [0,+\infty ]$ and $\psi :[0,+\infty ]\rightarrow [0,+\infty ]$
 be increasing functions for which there are $s_0, \,s_1\in {\mathbb{N}}$,  $s_0<s_1$, such  that 
\begin{itemize}
\item[{(i)}] the inverse functions $\varphi ^{-1}$ and  $\psi^{-1}$ are defined on $[\min \{\,\varphi (s_1-s_0),\psi (s_1-s_0)\,\},+\infty ]$, 
\item[{(ii)}] $2K(\psi ^{-1}\circ \varphi )(s-s_0)\leq (\psi ^{-1}\circ \varphi )(s)$ for all $s\geq s_1$,
\item[{(iii)}]
the function 
\[ [s_1+1,+\infty )\ni s\mapsto \frac{(\psi ^{-1}\circ \varphi )(s+1)}{(\psi ^{-1}\circ \varphi )(s)}\in {\mathbb{R}}\]
is bounded, 
\item[{(iv)}] the following integral is convergent:
\[ \int\limits_{s_1}^{+\infty }\frac{ds}{\varphi (s-s_0)^{1/(n-1)}}<+\infty .\]
\end{itemize}
Let ${\mathcal{F}}_K$ be a family of $K$-quasinearly subharmonic functions $u:\, D\rightarrow [-\infty ,+\infty )$ such that
\[u(x)\leq F_K(x)\]
for all  $x\in D$, where $F_K:\, D\rightarrow [0,+\infty ]$ is a Lebesgue measurable function. If for each compact set $E\subset D$,
\[\int\limits_E\psi (F_K(x))\, dm_n(x)<+\infty ,\]
then the family ${\mathcal{F}}_K$ is locally (uniformly) bounded in $D$. Moreover, the function $w^*:\, D\rightarrow [0,+\infty )$ is a $K$-quasinearly subharmonic function. Here
\[w^*(x):=\limsup_{y\rightarrow x}w(y),\]
where \[w(x):=\sup_{u\in {\mathcal{F}}_K}u^+(x).\]
}
\subsection{} The proof  of the theorem will be based on the following lemma, which has its origin in \cite{Do57},  Lemma~1, pp.~431-432, see also \cite{ArGa93}, Proposition~2, pp.~257-259. Observe that  we have applied our rather general and flexible lemma already  before (unlike previously, now we allow also the value $+\infty$ for our ``test functions'' $\varphi$ and $\psi$; observe that this does not cause any changes in the proof of our lemma, see \cite{Ri08$_1$}, pp.~5-8) when considering quasinearly subharmonicity of separately quasinearly subharmonic functions. As a matter of fact,  this lemma enabled us to  slightly improve Armitage's and Gardiner's almost sharp condition, see \cite{ArGa93}, Theorem~1, p.~256,  which ensures a separately subharmonic function to be subharmonic. See  \cite{Ri08$_1$}, Corollary~4.5, p.~13, and \cite{Ri08$_2$}, Corollary~3.3.3, p.~2622.   
\subsection*{Lemma} (\cite{Ri08$_1$}, Lemma~3.2, p.~5, and Remark~3.3, p.~8) \emph{Let $K$, $\varphi$, $\psi$ and $s_0, s_1\in {\mathbb{N}}$ be as in Theorem~2.1.
Let $u:\, D\rightarrow [0,+\infty )$ be a $K$-quasinearly subharmonic function. Let $\tilde {s}_1\in {\mathbb{N}}$, $\tilde {s}_1\geq s_3$, be arbitrary, where $s_3:=\max\{\, s_1+3,(\psi ^{-1}\circ \varphi )(s_1+3)\,\}$. 
Then for each $x\in D$ and $r>0$ such that 
$\overline{B^n(x,r)}\subset D$ either
\[u(x)\leq (\psi ^{-1}\circ \varphi )(\tilde {s}_1+1)\]
or
\[\Phi (u(x))\leq \frac{C}{r^n}\int\limits_{B^n(x,r)}\psi (u(y))\, dm_n(y)\]
where $C=C(n,K,s_0)$ and 
$\Phi :\, [0,+\infty )\rightarrow 
[0,+\infty ),$
\begin{displaymath}
 \Phi (t):=\begin{cases} \left[\int\limits_{(\varphi  ^{-1}\circ \psi )(t)-2}^{+\infty }\frac{ds}{\varphi (s-s_0)^{1/(n-1)}}\right]^{1-n},& 
{\textrm{when }} t\geq s_3, \\
\frac{t}{s_3}\Phi (s_3), &{\textrm{when }} 0\leq t<s_3.
\end{cases}
\end{displaymath}
}

\vspace{1ex}

\subsection{Proof of Theorem~2.1} 
Let $E$ be an arbitrary compact subset of $D$. Write $\rho_0:=$dist$(E,\partial D)$. Clearly $\rho_0>0$. Write 
\begin{displaymath}E_1:=\bigcup_{x\in E}\overline{B^n(x,\frac{\rho_0}{2})}.\end{displaymath}
Then $E_1$ is compact, and $E\subset E_1\subset D$.  Take $u\in {\mathcal{F}}_K^+$ arbitrarily,
where
\[{\mathcal{F}}_K^+:=\{\,u^+\, :\, u\in {\mathcal{F}}_K\, \}.\]
Let $\tilde{s}_1=s_1+2$, say. Take $x\in E$ arbitrarily and suppose that $u(x)> \tilde{s}_3$, where $\tilde{s}_3:=\max \{\, \tilde{s}_1+3,(\psi^{-1}\circ \varphi )(\tilde{s}_1+3)\,\}$, say. Using our Lemma and the assumption, we get
\begin{align*}
\left(\int\limits_{(\varphi^{-1}\circ \psi )(u(x))-2}^{+\infty}\frac{ds}{\varphi (s-s_0)^{1/(n-1)}}\right)^{1-n}&\leq \frac{C}{\left(\frac{\rho_0}{2}\right)^n}\int\limits_{B^n(x,\frac{\rho_0}{2})}\psi (u(y))\, dm_n(y)\\
&\leq \frac{C}{\left(\frac{\rho_0}{2}\right)^n}\int\limits_{E_1}\psi (F_K(y))\, dm_n(y)
<+\infty .
\end{align*}
Since
\begin{equation*}
\int\limits_{s_1}^{+\infty}\frac{ds}{\varphi (s-s_0)^{1/(n-1)}} 
<+\infty 
\end{equation*}
and $1-n<0$, the set of values 
\begin{equation*}
(\varphi^{-1}\circ \psi )(u(x))-2, \quad x\in E, \quad u\in {\mathcal{F}}_K^+,
\end{equation*}
is  bounded. Thus also the set of values
\begin{equation*}
u(x), \quad x\in E, \quad u\in {\mathcal{F}}_K^+,
\end{equation*}
is bounded.

To show that $w^*$ is $K$-quasinearly subharmonic in $D$,  proceed as follows. Take $x\in D$ and $r>0$ such that $\overline{B^n(x,r)}\subset D$. 
For each $u\in {\mathcal{F}}^+_K$ we have then
\begin{equation*}
u(x)\leq \frac{K}{\nu_n r^n}\int\limits_{B^n(x,r)}u(y)\, dm_n(y).\end{equation*}  
Since 
\begin{equation*}
u(x)\leq \sup_{u\in {\mathcal{F}}_K^+}u(x)=w(x)\leq w^*(x),\end{equation*}
we  have 
\begin{equation}
w(x)\leq \frac{K}{\nu_n r^n}\int\limits_{B^n(x,r)}w^*(y)\,dm_n(y).
\end{equation}

Then just take the upper semicontinuous regularizations on  both sides of (5) and use Fatou Lemma on the right hand side (this is of course possible, since  $w^*$ is locally bounded in $D$), say: 
\begin{align*}
\limsup_{y\rightarrow x}w(y)&\leq \limsup_{y\rightarrow x}\frac{K}{\nu_n r^n}\int\limits_{B^n(y,r)}w^*(z)\, dm_n(z)\\
&\leq \limsup_{y\rightarrow x}\frac{K}{\nu_n r^n}\int w^*(z)\,\chi_{\overline{B^n(y,r)}}(z)\, dm_n(z)\\
&\leq \frac{K}{\nu_n r^n}\int w^*(z)\left( \limsup_{y\rightarrow x}\chi_{\overline{B^n(y,r)}}(z)\right) \, dm_n(z).
\end{align*}
Since for all $z\in D$,
\begin{equation*}
\limsup_{y\rightarrow x}\chi_{\overline{B^n(y,r)}}(z)\leq \chi_{\overline{B^n(x,r)}}(z),
\end{equation*}
we get the desired inequality
\begin{equation*}
w^*(x)\leq \frac{K}{\nu_n r^n}\int\limits_{B^n(x,r)}w^*(y)\, dm_n(y).
\end{equation*} 
\null{} \hfill $\qed$
\subsection{Remark} If $w$ is Lebesgue measurable, it follows that already $w$ is $K$-quasinearly subharmonic.
\subsection{Corollary}  \emph{Let $\varphi :[0,+\infty ]\rightarrow [0,+\infty ]$ be a strictly increasing function such that for some $s_0, s_1\in {\mathbb{N}}$, $s_0<s_1$,
\begin{equation*}\int\limits_{s_1}^{+\infty} \frac{ds}{[\varphi (s-s_0)]^{1/(n-1)}}<+\infty .
\end{equation*}
Let ${\mathcal{F}}_K$ be a family of $K$-quasinearly subharmonic functions $u:\, D\rightarrow [-\infty ,+\infty )$ such that
\[u(x)\leq F_K(x)\]
for all  $x\in D$, where $F_K:\, D\rightarrow [0,+\infty ]$ is a Lebesgue measurable function. Let  $p>0$ be arbitrary. If for each compact set $E\subset D$,
\begin{equation*}\int\limits_E\varphi ( \log^+[F(x)]^p)\, dm_n(x)<+\infty ,
\end{equation*}
then the family ${\mathcal{F}}_K$ is locally (uniformly) bounded in $D$. Moreover, the function $w^*:\, D\rightarrow [0,+\infty )$ is a $K$-quasinearly subharmonic function. Here
\[w^*(x):=\limsup_{y\rightarrow x}w(y)\]
where \[w(x):=\sup_{u\in {\mathcal{F}}_K}u^+(x).\]
}
For the proof, take $p>0$ arbitrarily,  choose   $\psi (t)=(\varphi \circ \log^+)(t^p)$, and just check  that the conditions (i)--(iv) indeed hold. 

The case $p=1$ and $K=1$ gives  Domar's and Rippon's results, Theorems~A and B above.  
\subsection{Remark} As a matter of fact, our result includes even the case when $\psi (t)=(\varphi \circ \log^+)(\phi(t))$, where $\phi :\, [0,+\infty ]\rightarrow [0,+\infty ]$ is \emph{any} strictly increasing function which satisfies the following two conditions:
\begin{itemize}
\item[(a)] $\phi^{-1}$ satisfies the $\Delta_2$-condition, 
\item[(b)] $2K \phi^{-1}(e^{s-s_0})\leq \phi^{-1}(e^s)$ for all $s\geq s_1$.
\end{itemize}  

\end{document}